\documentclass{amsart}
\usepackage{amsmath,latexsym}
\usepackage[mathscr]{eucal}
\usepackage{amssymb}

\newtheorem{theorem}{Theorem}%[section]

\numberwithin{equation}{section}
 \allowdisplaybreaks \theoremstyle{remark}

\def\star{\raise .5ex \hbox{*}}
\def\sumstar_#1{\setbox0=\hbox{$\scriptstyle{#1}$}
\setbox2=\hbox{$\displaystyle{\sum}$}
\setbox4=\hbox{${}\star\mathsurround=0pt$} \dimen0=.5\wd0
\advance\dimen0 by-.5\wd2 \ifdim\dimen0>0pt \ifdim\dimen0>\wd4
\kern\wd4 \else\kern\dimen0\fi\fi \mathop{{\sum}\star}_{\kern-\wd4
#1}}

\begin{document}

\title{Positive Proportion of Small Gaps Between Consecutive Primes}
\author{D. A. Goldston, J. Pintz and C. Y. Y{\i}ld{\i}r{\i}m }

\thanks{The first author was supported by NSF; this study was done while
the second and third authors were members at the Institute for
Advanced Studies, Princeton during Fall 2009 and were supported by
the Oswald Veblen Fund; the second author also acknowledges the
partial support of ERC-AdG. No.228005.}

\date{\today}

\maketitle

\vskip .3in

\section{Introduction } Let $\pi(x)$ denote as usual the number of primes $\le x$. The prime number theorem is the asymptotic relation $\pi(x) \sim \frac x{\log x}$ as $x\to \infty$. Now let $p_{n}$ be the $n^{\text{th}}$ prime. We consider the gaps $p_{n+1}-p_n$ in the sequence of primes. By the prime number theorem the average of this sequence of gaps is $\log p_n$. In \cite{GPY1} we proved that
\begin{equation}\label{eq:1.1}
\liminf_{n\to\infty}{p_{n+1}-p_{n}\over \log p_{n}} =0, 
\end{equation}
so that there are gaps arbitrarily smaller than the average.
In this paper we prove that these small gaps occur so frequently that they form a positive proportion of all the gaps.  Define the distribution function for small gaps between consecutive primes by
\begin{equation} \label{eq:1.2}  
P(x, \eta) := \frac{1}{\pi(x)} \sum_{ \substack{p_n\le x\\ p_{n+1}-p_n \le \eta \log p_n}}1. \end{equation}

\begin{theorem} For any fixed $\eta >0$, we
have 
\begin{equation} \label{eq:1.3}
P(x,\eta)  \gg_{\eta} 1, \qquad \text{as} \ \  x\to \infty. 
\end{equation}
Thus the small gaps between consecutive primes constitute a positive proportion of the set of all
gaps between consecutive primes. 
\end{theorem}

Our method actually obtains an explicit dependence on $\eta$ in the lower bound in \eqref{eq:1.3}; we leave this for a later paper \cite{GPY4}. (The result is exponentially small in a power of $\eta$.) It has been conjectured that primes are distributed around their average spacing in a Poisson distribution. Gallagher \cite{Ga} has proved that this is a consequence of the Hardy-Littlewood prime tuple conjecture. If this is the case, then for fixed $\eta>0$,
\begin{equation} \label{eq:1.4} P(x,\eta) \sim 1 -e^{-\eta} \qquad \text{as} \ \  x \to \infty, \end{equation}
and consequently 
\begin{equation} \label{eq:1.5} P(x,\eta) \sim \eta  \qquad \text{when} \ \ \eta \to 0 \quad \text{sufficiently slowly as} \ \ x \to \infty. \end{equation}

By sieve methods it is easy to obtain an upper bound of this magnitude for $P(x,\eta)$. Obviously $P(x,\eta) \le 1$, so we only need to consider $\eta \ll 1$.  

\begin{theorem} For  $1/\log x \ll \eta \ll 1$, we
have 
\begin{equation} \label{eq:1.6}
P(x,\eta)  \ll  \eta, \qquad \text{as} \ \  x\to \infty. 
\end{equation}
\end{theorem}
We thus see that Theorem 1 can not continue to hold if $\eta \to 0$; there are not a positive proportion of prime gaps smaller than $\eta \log p_n$ if $\eta\to 0$ as $p_n\to \infty$.

The method we use for our results on small gaps between primes uses information on the distribution of primes in arithmetic progressions, specifically what is called an admissible level of distribution $\vartheta$. The precise definition may be found in \cite{GPY1}, but roughly this means that the primes less than or equal to $x$  are distributed evenly among the arithmetic progressions $a (\text{mod}\ q)$, $(a,q)=1$, for almost all progressions with $1\le q \le Q= x^{\vartheta -\epsilon}$. The unconditional results use $\vartheta =1/2$ which is known to hold by the Bombieri-Vinogradov theorem. If we assume some $\vartheta > 1/2$ holds then we obtain bounded gaps between primes. The method we use here does not lead to the conjectured number of such bounded gaps, but when applied to $P(x,\eta)$ does obtain lower bounds closer to the conjectured asymptotic relation in \eqref{eq:1.5}. 

\begin{theorem} Suppose some $\vartheta \ge \vartheta_0>\frac12$ is an admissible level of distribution for primes. Then there exists an integer $m(\vartheta_0)$ such that for  $1/\log x \ll \eta \ll 1$, we
have 
\begin{equation} \label{eq:1.7}
P(x,\eta)  \gg  \eta^{m(\vartheta_0)}, \qquad \text{as} \ \  x\to \infty. 
\end{equation}
In particular if $\vartheta_0 > .971$ then $P(x,\eta)  \gg  \eta^{5}$ and if $\vartheta_0 > .953$ then $P(x,\eta)  \gg  \eta^{6}$. 
\end{theorem}

 Our method only applies to obtain results on pairs of nearby primes unless we assume the  Elliott-Halberstam conjecture that $\vartheta=1$ is admissible. With this conjecture the method is able to prove the existence of triples of primes closer than any multiple of the average spacing, although it can not produce bounded gaps between such triples. With this conjecture, the proof of Theorem 1 immediately leads to the following result.

\begin{theorem} Assume the Elliott-Halberstam conjecture that $\vartheta =1$ is an admissible level of distribution for primes. Then for any fixed $\eta>0$, we
have 
\begin{equation} \label{eq:1.8}
P_2(x,\eta) := \frac{1}{\pi(x)} \sum_{ \substack{p_n\le x\\ p_{n+2}-p_n \le  \eta \log p_n}}1 \gg_\eta 1, \qquad \text{as} \ \  x\to \infty. 
\end{equation}
\end{theorem}

\section{ Results from the paper Primes in Tuples I}

Let us recall briefly how the result \eqref{eq:1.1} was obtained. Consider
the $k$-tuple
\begin{equation}\label{eq:2.1}
\mathcal H = \{h_{1}, h_{2}, \ldots, h_{k}\} \;\; {\rm with \,\,
distinct \,\, integers} \;\; 1\le h_{1}, \ldots, h_{k} \le h  ,
\end{equation}
and for a prime $p$ denote by $\nu_{p}(\mathcal H)$ the number of
distinct residue classes modulo $p$ occupied by the entries of
$\mathcal H$. The singular series associated with $\mathcal H$ is
defined as
\begin{equation} \label{eq:2.2}
\mathfrak S (\mathcal H) : = \prod_{p}\left(1-{1\over
p}\right)^{-k}\left(1-{\nu_{p}(\mathcal H)\over p}\right) , 
\end{equation}
the product being convergent because $\nu_{p}(\mathcal H)=k$ for
$p>h$. We say that $\mathcal H$ is admissible if
\begin{equation} \label{eq:2.3}
P_{\mathcal H}(n) := (n+h_{1})(n+h_{2})\cdots (n+h_{k}) 
\end{equation}
is not divisible by a fixed prime number for every $n$, which is
equivalent to $\nu_{p}(\mathcal H) \neq p$ for all $p$ and therefore
also to $\mathfrak S(\mathcal H) \neq 0$. That $\{ n+h_1 , n+h_2 ,
\ldots, n+h_k \}$ is a prime tuple, i.e. each entry is prime, is
equivalent to $P_{\mathcal H}(n)$ being a product of $k$ primes.

Next, define the function 
\begin{equation} \label{eq:2.4}
\Lambda_{R}(n; \mathcal H, \ell) := {1\over
(k+\ell)!}\sum_{\substack{d\mid P_{\mathcal H}(n) \\ d \leq
R}}\mu(d)(\log {R\over d})^{k+\ell} , \quad ( k= |\mathcal H |, \ \ 0 \leq \ell < k) ,
\end{equation}
which is designed to approximate the indicator function for when $P_{\mathcal H}(n)$ has at most $k+\ell$ distinct prime factors. Let
\begin{equation} \label{eq:2.5}
\theta(n) := \begin{cases}
\log n &\text{ if $n$ is prime}, \\
0 &\text{ otherwise}.
\end{cases}
\end{equation} 
Now, as a consequence of Propositions 1 and 2 in \cite{GPY1}, we have the following three results. 
For $\mathcal H$ admissible (or equivalently $\mathfrak{S}(\mathcal{H}) \neq 0$), $h\le R \ll N^{\frac{1}{2}} (\log N)^{-B(k)}$ and $R,N\to \infty$, we have 
\begin{equation}\label{eq:2.6}
\sum_{n\leq N} \Lambda_{R} (n; \mathcal H, \ell)^2 \sim
\frac1{(k + 2\ell)!} {2\ell\choose \ell} \mathfrak S(\mathcal
H) N (\log R)^{k + 2\ell}.
\end{equation}
For any  $h_i \in \mathcal H$ and $\mathfrak{S}(\mathcal{H}) \neq 0$, we have for $h\le R^\epsilon$, $R\ll N^{\frac{\vartheta}{2} -\epsilon}$, and $R,N\to \infty$, 
\begin{equation} \label{eq:2.7}
\sum_{n\leq N} \Lambda_R(n; \mathcal H , \ell)^2 \theta (n +
h_i) \sim \frac1{(k + 2\ell + 1)!} {2\ell + 2\choose \ell+1}
\mathfrak S(\mathcal H) N(\log R)^{k+2\ell+1},
\end{equation}
and for $h_0\not \in \mathcal{H}$ and $\mathfrak{S}(\mathcal{H}\cup \{h_0\})) \neq 0$,
\begin{equation} \label{eq:2.8}
\sum_{n\leq N} \Lambda_R(n; \mathcal H, \ell)^2\,  \theta (n +
h_0) \sim \frac1{(k + 2\ell )!} {2\ell \choose \ell}
\mathfrak S(\mathcal{H}\cup \{h_0\}) N(\log R)^{k+2\ell}.
\end{equation} 

We also need a result of Gallagher \cite{Ga}: as $h\to \infty$,
\[
\sum_{\substack{1\le h_1,h_2, \ldots ,h_k\le h \\
\text{distinct}}}
\mathfrak{S}(\mathcal{H} ) \sim h^k.
\]
However, we  now change notation slightly from \cite{GPY1}. Equation \eqref{eq:2.1} is equivalent to the conditions that $|\mathcal H|=k$ and $\mathcal H \subset \{ 1, 2, \ldots , \lfloor h\rfloor\}$. Further, Gallagher's result is unchanged if we restrict ourselves to the non-zero terms where $\mathcal H$ is admissible. Hence Gallagher's result can be restated as
\begin{equation} \label{eq:2.9}
\sum_{\substack{ |\mathcal H|=k\\ \mathcal H \subset \{1, 2, \ldots , \lfloor h\rfloor \} \\ \mathcal H \ \text{admissible}}} \mathfrak{S}(\mathcal{H} ) \sim \frac{h^k}{k!}, 
\end{equation}
where the $k!$ is from the permutation of the elements of $\mathcal H$ which we no longer sum over. 
Now, define, for $\nu$ a positive integer which in this paper is either 1 or 2, 
\begin{equation} \label{eq:2.10}
\mathcal{S}
:= \sum_{\substack{ |\mathcal H|=k\\ \mathcal H \subset \{1, 2, \ldots , \lfloor h\rfloor \} \\ \mathcal H \ \text{admissible}}} \left(\sum_{n=N+1}^{2N}
\Bigg( \sum_{\substack{1\le h_0\le h \\  \mathfrak{S}(\mathcal H \cup \{h_0\}) \neq 0}} \theta (n+h_0) - \nu \log 3N \Bigg)
\Lambda_{R}(n;\mathcal{H}, \ell)^2 \right).
\end{equation}
Applying \eqref{eq:2.6}--\eqref{eq:2.9} (and noting that once these equations are used the conditions on admissibility may be dropped),  a simple calculation gives 
\begin{equation}\begin{split} \label{eq:2.11}
 \mathcal{S}&\sim \sum_{\substack{ |\mathcal H|=k\\ \mathcal H \subset \{1, 2, \ldots , \lfloor h\rfloor \} }} 
\Bigg(\frac{k}{(k+2\ell+1)!}\genfrac{(}{)}{0pt 
}{0}{2\ell+2}{\ell+1} \mathfrak{S}(\mathcal{H}) N(\log 
R)^{k+2\ell+1} 
\\ 
& \quad +\sum_{\substack{1\le h_0\le h \\ h_0\neq h_i, 1\le i\le k 
}} \frac{1}{(k+2\ell)!}\genfrac{(}{)}{0pt }{0}{2\ell}{\ell} 
\mathfrak{S}(\mathcal{H}\cup \{h_0\}) N(\log R)^{k+2\ell} \\ 
& \qquad  - \nu \log 3N  \frac{1}{(k+2\ell)!}\genfrac{(}{)}{0pt 
}{0}{2\ell}{\ell} \mathfrak{S}(\mathcal{H}) N(\log 
R)^{k+2\ell}\Bigg)\\ 
& \sim  \mathcal M(k,\ell, h) \frac{1}{(k+2\ell)!k!}\genfrac{(}{)}{0pt 
}{0}{2\ell}{\ell} Nh^k(\log R)^{k+2\ell},
\end{split}\end{equation}
where 
\begin{equation} \label{eq:2.12} 
\mathcal M(k,\ell, h) := \frac{2k}{k+2\ell+1}\frac{2\ell+1}{\ell +1} 
\log R + h - \nu \log 
3N. 
\end{equation}
(Note that in the calculation above each of the sets $\mathcal H \cup \{ h_0\}$ occurred $k+1$ times in the summation.)
 Thus, there
are at least $\nu +1$ primes in some interval $(n,n+h]$, $N<n\le 2N$, provided that $\mathcal M(k,\ell, h)>0$. Taking $R= N^{\frac{\vartheta}{2} - \epsilon}$, this is true when
\begin{equation}\label{eq:2.13}
h  > \biggl(\nu -\frac{2k}{k+2\ell+1}\frac{2\ell+1}{\ell +1}
\biggl(\frac{\vartheta}{2}-\epsilon\biggr) \biggr)\log 3N,
\end{equation}
which, on letting $\ell = \lfloor \frac{\sqrt{k}}2 \rfloor$ and taking $k$ sufficiently large,
gives
\begin{equation}\label{eq:2.14} h> \left(\nu -  2\vartheta  +4 \epsilon +O\biggl(\frac{1}{\sqrt{k}}\biggr) \right)\log N.
\end{equation}
Taking $\nu =1$ and $\vartheta=1/2$ proves \eqref{eq:1.1}. 

\section{A new prime tuple detecting weight.}
The prime pairs we found in the last section are counted with the weight
$\Lambda_{R}(n;\mathcal{H}, \ell)^2$,
and this weight needs to be removed in order to count the number of prime pairs themselves.
As usual, this is accomplished by using Cauchy's inequality. The problem with this approach (which stumped 
us for many years) is that there are values of $n$ with many divisors for which $\Lambda_{R}(n;\mathcal{H}_k, \ell)^2$ 
is exceptionally large, and these terms prevent us from obtaining the desired positive proportion result. The solution of this problem
was found by Pintz in \cite{P}, and is based on a general property of the Selberg sieve. This property is that the Selberg sieve weights effectively remove almost
all numbers with many prime factors. Therefore the $n$ for which $\Lambda_{R}(n;\mathcal{H}_k, \ell)^2$ may be large are also numbers which contribute very little to the total size of the asymptotic formulas in \eqref{eq:2.6}--\eqref{eq:2.8}.

 We define
\begin{equation} \label{eq:3.1}
\mathcal{P}(x) := \prod_{p_{n}\leq x}p_n .
\end{equation}
Let $\delta >0$ be a fixed constant that we can choose to be as small as we wish. We want to remove from our earlier sums  the terms when $(P_{\mathcal H} (n),\mathcal P(R^\delta)) >1$.   We can do this with an error that is small when $\delta$ is small by Pintz's work \cite{P}. The results we need are immediate consequences of Pintz's Lemmas 4 and 5 and \eqref{eq:2.6}--\eqref{eq:2.8}.
We take $\ell \asymp \sqrt{k}$ which eliminates the $\ell$ dependence in the error terms which follow.
Suppose  $N^{c_1} \leq R \leq
N^{{1\over 2+\delta}}(\log N)^{-C_1}$ where $c_1$ and $C_1$ are
suitably chosen constants depending on $k$. (Actually $c_1={1\over 5}$ and $C_1$ taken sufficiently large suffices.)  If $\mathcal H$ is admissible with $h \ll \log R$ and
$h\to\infty$ with $N$, we have
\begin{equation} \label{eq:3.2}
\sum_{\substack{n=N+1 \\ (P_{\mathcal 
H}(n), \mathcal{P}(R^{\delta}))>1} }^{2N}\Lambda_R( n;\mathcal{H},\ell)^2
\ll_{k} \delta\  \mathfrak S(\mathcal H) N (\log R)^{k+2\ell }.
\end{equation}
For $1\leq h_0 \leq h$, write $m=1$ when $h_0 \in \mathcal
H$ and $m=0$ when $h_0 \not\in \mathcal H$. For $\epsilon >0$ if  $\mathfrak{S}(\mathcal H \cup \{h_0\})\neq 0$, then for $N^{c_1} \leq R \leq N^{\frac{\vartheta - \epsilon}{(2+\delta)}}$  we have
\begin{equation} \label{eq:3.3}
 \sum_{\substack{n=N+1 \\ (P_{\mathcal 
H}(n), \mathcal{P}(R^{\delta}))>1} }^{2N}\theta
(n+h_0)\Lambda_R (n;\mathcal{H},\ell)^2 \ll_k \delta \ \mathfrak S(\mathcal H\cup \{h_0\}) N(\log R)^{k + 2\ell +m}.
\end{equation}
We now define a modified prime tuple approximation weight 
\begin{equation} \label{eq:3.4}
{\Lambda^{\!*}}_{\!\!\!R}( n;\mathcal{H},\ell,\delta) := \begin{cases}
 \Lambda_R( n;\mathcal{H},\ell)&\text{ if $(P_{\mathcal 
H}(n), \mathcal{P}(R^{\delta}))=1$}, \\
0 &\text{ otherwise}.
\end{cases}
\end{equation}
We thus see that this weight tries to approximate prime tuples using only almost prime divisors, and is only insignificantly less effective than our original approximation $\Lambda_R$ when $\delta $ is taken sufficiently small. 

\section{ Detecting Pairs of Primes using the new approximation}
We now replace $\mathcal S$ in \eqref{eq:2.10} with
\begin{equation} \label{eq:4.1}
\mathcal{S}^*
:= \sum_{\substack{ |\mathcal H|=k\\ \mathcal H \subset \{1, 2, \ldots , \lfloor h\rfloor \} \\ \mathcal H \ \text{admissible}}} \left(\sum_{n=N+1}^{2N}
\Bigg( \sum_{\substack{1\le h_0\le h \\  \mathfrak{S}(\mathcal H \cup \{h_0\}) \neq 0}} \theta (n+h_0) - \nu \log 3N \Bigg)
{\Lambda^{\!*}}_{\!\!\!R}( n;\mathcal{H},\ell,\delta)^2 \right).
\end{equation}
Using \eqref{eq:3.3} and \eqref{eq:3.4}, the difference $\mathcal S^*-\mathcal S$ is  
\[\begin{split} 
  &\ll_k \delta \sum_{\substack{ |\mathcal H|=k\\ \mathcal H \subset \{1, 2, \ldots , \lfloor h\rfloor \} }} 
\bigg( \mathfrak{S}(\mathcal{H})\big(\log R +\log 3N\big)   +\sum_{\substack{1\le h_0\le h \\ h_0\neq h_i, 1\le i\le k 
}} \mathfrak{S}(\mathcal{H}\cup \{h_0\}) \bigg)N(\log 
R)^{k+2\ell}\\ 
& \ll_k \delta Nh^k(\log R)^{k+2\ell}\log N,
\end{split}\]
where we used that $h\ll \log R$. 
We conclude from this and \eqref{eq:2.11} that
\begin{equation} \label{eq:4.2}
 \mathcal{S}^*\sim  \left(\mathcal M(k,\ell, h) +O_k(\delta \log N) \right)\frac{1}{(k+2\ell)!k!}\genfrac{(}{)}{0pt 
}{0}{2\ell}{\ell} Nh^k(\log R)^{k+2\ell},
\end{equation}
as $R,N\to \infty$, where  $N^{c_1} \leq R \leq N^{\frac{\vartheta - \epsilon}{(2+\delta)}}$.  Clearly if $h= \eta \log N$ and  $R$, $\vartheta$, and $k$ are chosen appropriately to make $\mathcal M(k,\ell,h) $ positive we can then choose $\delta$ sufficiently small  so that $\mathcal{S}^*$ will also be positive. Then as in Section 2 we will have produced pairs of nearby primes. 

\section{ Removing the weight}

The property that ${\Lambda^{\!*}}_{\!\!\!R}( n;\mathcal{H},\ell,\delta)$ possesses that $\Lambda_R(n;\mathcal H, \ell)$ lacks is that  it is never larger than some constant depending on $k$ and $\delta$ times the size of the single term in its sum from the divisor $d=1$. To see this, note that  all prime
factors of $P_{\mathcal H}(n)$ in the sum that forms ${\Lambda^{\!*}}_{\!\!\!R}( n;\mathcal{H},\ell,\delta)$ are greater than $R^{\delta}$, and thus the number of squarefree divisors of
$P_{\mathcal H}(n)$ is at most $2^{{k\log 3N\over \delta\log R}}$. Thus for $N^{c_1}\le R$, 
\begin{equation} \label{eq:5.1}
{\Lambda^{\!*}}_{\!\!\!R}( n;\mathcal{H},\ell,\delta) \leq {2^{{k\log 3N\over \delta\log
R}}\over (k+\ell)!}(\log R)^{k+\ell} \ \ll_{k,\delta}  (\log R)^{k+\ell}.
\end{equation}

We now proceed to obtain an upper bound for $\mathcal{S}^*$ which counts small gaps between consecutive primes without weights.  First, letting 
\begin{equation} \label{eq:5.2} 
\Theta(n,h) := \sum_{1\le h_0\le h}\theta(n+h_0), \quad  \pi(n,h) := \pi(n+h)-\pi(n),
\end{equation}
we have
\begin{equation} \label{eq:5.3}
\begin{split}
\mathcal{S}^*
&\leq  \sum_{\substack{ |\mathcal H|=k\\ \mathcal H \subset \{1, 2, \ldots , \lfloor h\rfloor \} \\ \mathcal H \ \text{admissible}}} \Bigg(\sum_{\substack{n=N+1\\ \Theta(n,h) > \nu \log 3N}}^{2N}
\Theta(n,h)
{\Lambda^{\!*}}_{\!\!\!R}( n;\mathcal{H},\ell,\delta)^2 \Bigg)\\
& \ll_{k,\delta} (\log R)^{2k+2\ell} \log 3N \sum_{\substack{n=N+1\\ \pi(n,h) > \nu }}^{2N} \pi(n,h) \sum_{\substack{ |\mathcal H|=k \\ \mathcal H \subset \{1, 2, \ldots , \lfloor h\rfloor \} \\ \mathcal H \ \text{admissible}\\ (P_{\mathcal H}(n), \mathcal P(R^\delta))=1}} 1  .
\end{split}
\end{equation}
Denote the inner sum by $ \mathcal{T}(\mathcal{H},n)$, and let
\begin{equation} \label{eq:5.4}
Q_\nu(N,h):= \sum_{\substack{n= N\\
\pi(n,h)> \nu} }^{2N}1 .
\end{equation}
We now have by Cauchy's inequality that 
\begin{equation} \label{eq:5.5}
 \sum_{\substack{n=N+1\\ \pi(n,h) > \nu}}^{2N} \pi(n,h)\mathcal{T}(\mathcal{H},n) \le  Q_\nu(N,h)^{\frac12}\left(\sum_{n=N+1}^{2N} \pi(n,h)^2 \mathcal{T}(\mathcal{H},n)^2\right)^{\frac12}.  
\end{equation}
If $n$ is an integer for which $\pi(n+h)-\pi(n) > \nu$, then there
must be a $j$ such that $n < p_j$ and $p_{j+\nu} \leq n+h$. Thus
$p_{j+\nu}-p_j < h$ and $p_{j+\nu} -h \leq n < p_j $, so that
there are   less than $\lfloor h \rfloor $ such integers $n$
corresponding to each such gap. Therefore
\begin{equation} \label{eq:5.6}
Q_\nu(N,h) \leq \, h\sum_{\substack{ N < p_{j}\le 2N \\
p_{j+\nu}-p_j \leq h}}1 \, + O(Ne^{-c\sqrt{\log N}}),
\end{equation}
where we have used the prime number theorem with error term to
remove the prime gaps which overlap the endpoints. (This is
explicitly shown in \cite{GY1}). We will prove below that, for $2\le h\le \log R$,
\begin{equation} \label{eq:5.7}
\sum_{n=N+1}^{2N} \pi(n,h)^2 \mathcal{T}(\mathcal{H},n)^2 \ll_{k,\delta} \left(\frac{h}{\log R}\right)^k N,
\end{equation}
which on combining with \eqref{eq:5.3} and \eqref{eq:5.5} produces the upper bound 
\begin{equation} \label{eq:5.8}
\mathcal{S}^* \ll_{k,\delta}  (\log R)^{2k+2\ell} (\log 3N)  Q_\nu(N,h)^{\frac12} \left(\frac{h}{\log R}\right)^{\frac{k}{2}}N^{\frac12} 
\end{equation}
Together with \eqref{eq:5.6} this provides the desired lower bound for the unweighted number of small prime gaps. 

To prove \eqref{eq:5.7}, we recall that the main theorem of Selberg's upper bound sieve
(Theorem 5.1 of \cite{HR} or Theorem 2 in \S 2.2.2 of \cite{Gr})
gives for any set $\mathcal H$ and $\delta < \frac12$
\begin{equation} \label{eq:5.9}
\sum_{\substack{ n=N+1\\
(P_{\mathcal{H}}(n), \mathcal{P}(R^{\delta}))=1}}^{2N}1 \leq
\frac{|\mathcal{H}|!\mathfrak S (\mathcal{H})}{(\log
R^{\delta})^{|\mathcal{H}|}}N(1+o(1)), \quad \; (N\to\infty).
\end{equation}
Writing 
\[ \pi(n,h) = \sum_{\substack{1\le h'\le h\\  n+h' \ \text{prime}}} 1,\]
we see that the left-hand side of \eqref{eq:5.7} is
\[ \ll \sum_{1\le h',h'' \le h} \sum_{\substack{ |\mathcal H_i|=k \\ \mathcal H_i \subset \{1, 2, \ldots , \lfloor h\rfloor \} \\ \mathcal H_i \ \text{admissible}\\ i=1,2}} \ \ \sum_{\substack{n=N+1 \\ (P_{\mathcal 
H_1}(n), \mathcal{P}(R^{\delta}))=1 \\ (P_{\mathcal 
H_2}(n), \mathcal{P}(R^{\delta}))=1\\ n+h' , \ n+h''\ \text{prime}} }^{2N}  1,
\]
The conditions on the inner sum are weakened if we let $\mathcal H_0 = \{h'\}\cup \{h''\} \cup \mathcal H_1 \cup \mathcal H_2$ and require $ (P_{\mathcal 
H_0}(n), \mathcal{P}(R^{\delta}))=1$, and therefore we obtain the upper bound
\[ \sum_{r=k}^{2k+2} \sum_{\substack{ |\mathcal H_0|=r \\ \mathcal H_0 \subset \{1, 2, \ldots , \lfloor h\rfloor \} }} \ \ \sum_{\substack{n=N+1 \\ (P_{\mathcal 
H_0}(n), \mathcal{P}(R^{\delta}))=1}}^{2N} 1.\]
By \eqref{eq:5.9} and \eqref{eq:2.9} this is, for $2\le h\le \log R$,
\[ \ll_{k} \sum_{r=k}^{2k+2} \sum_{\substack{ |\mathcal H_0|=r \\ \mathcal H_0 \subset \{1, 2, \ldots , \lfloor h\rfloor \} }} \frac{\mathfrak{S}(\mathcal H_0)}{(\log R^\delta)^r} N \quad \ll_{k,\delta} N \sum_{r=k}^{2k+2} \left(\frac{h}{\log R}\right)^r \quad \ll_{k,\delta} \left(\frac{h}{\log R}\right)^k N, \]
which is \eqref{eq:5.7}.

\section{Proof of the Theorems}

We now take $R= N^{\frac{\vartheta - \epsilon}{(2+\delta)}}$, $h= \eta \log N$, and $\frac2{\log N}\le \eta \le \frac15$ so that $h\le \log R$. Combining \eqref{eq:4.2} and \eqref{eq:5.8} we obtain
\begin{equation} \label{eq:6.1} 
\left(\mathcal M(k,\ell , h) + O_k(\delta \log N) \right) \left(\frac h{\log R}\right)^{\frac k2}\frac{N^{\frac12}}{\log N} \le C(k,\delta)  Q_\nu(N,h)^{\frac 12},
\end{equation}
where $C(k,\delta) >0$ is a (large) constant depending on $k$ and $\delta$.

We first prove Theorems 1 and 4. Taking $\ell = \lfloor \frac{\sqrt{k}}2 \rfloor$, we find
\[ \frac{2k}{k+2\ell+1}\frac{2\ell+1}{\ell +1} > 4 - \frac{c_2}{\sqrt{k}}, \quad k\ge 4, \]
for a suitable constant $c_2$ (A short calculation shows $c_2 =8$ works here.)
Hence from \eqref{eq:2.12}  we have
\begin{equation}\label{eq:6.2}
\begin{split} \mathcal{M}(k,\ell, h) + O_k(\delta \log N) &>  \left(4 - \frac{c_2}{\sqrt{k}}\right)\left(\frac{\vartheta -\epsilon}{2 + \delta}\right)\log N  + h -\nu \log 3N - c_3(k)\delta \log N \\&
> \left(\eta + (2\vartheta - \nu) - 4\epsilon - \frac{c_2}{\sqrt{k}}- c_4(k)\delta\right) \log N. \end{split}
\end{equation}
Take $\vartheta=\frac12$ and $\nu=1$ for Theorem 1, or $\vartheta=1$ and $\nu=2$ for Theorem 4. Hence, given a fixed $\eta >0$ we can first choose $k=k(\eta)$ large enough and then $\epsilon =\epsilon(\eta)$ and $\delta = \delta(\eta)$ small enough so that
\[ \mathcal{M}(k,\ell, h) + O_k(\delta \log N) > \frac\eta 2 \log N.\]
From \eqref{eq:6.1} we immediately obtain 
 $Q_\nu(N,h) \gg_\eta N$, and  \eqref{eq:5.6} completes the proof. 

For the proof of  Theorem 3 we take $\nu =1$,  and note that  if $\vartheta \ge \vartheta_0 > \frac 12$,  then we do not need $\eta$ to make the right-hand side  \eqref{eq:6.2} positive; we only need to make $k$ large enough and then $\delta$ small enough to accomplish this. Hence, with $k=k_0(\vartheta_0)$, $\ell  =\ell_0(\vartheta_0)$, and $\delta = \delta_0(\vartheta_0)$ we have
\[ \mathcal{M}(k,\ell, h) + O_k(\delta \log N) \gg_{k,\delta}  \log N.\]
From \eqref{eq:6.1} we then obtain
 $Q_1(N,h) \gg \eta^{k(\vartheta_0)} N$, and then from \eqref{eq:5.6} the first part of Theorem 3 follows with $m(\vartheta_0) = k(\vartheta_0)-1$.
Next, we take $k=7$ and $\ell=1$ in 
\eqref{eq:2.12} and obtain
\[ \mathcal M(7,1,h) = \frac{21}{10}\log R +h - \log 3N\]
and this is $\gg \log N$ independent of $h$ provided $\vartheta\ge \vartheta_0$ where
\[2\left(\frac{\vartheta_0 -\epsilon}{2 + \delta}\right) > \frac{20}{21} = 0.95238\ldots .\] 
 Hence on taking $\epsilon$ and $\delta$ sufficiently small depending on $\vartheta_0> .953$, then for $N\ge N_0(\vartheta)$ we have as above $Q_1(N,h) \gg \eta^7 N$, and hence by \eqref{eq:5.6} we conclude $P(2N,\eta)\gg \eta^6$. For the final part of Theorem 3 we wish to take $k=6$ in \eqref{eq:2.12} as above but this just fails to give a positive result. However by using a linear combination of $\Lambda_R$'s with $k=6$ and $\ell=0$ and $\ell=1$ we are able to obtain a positive result here provided $\vartheta_0 > .971$ as was done in \cite{GPY1}. The proof then follows as above with minor changes. 

Finally, we prove Theorem 2. This is an almost immediate consequence of any sieve upper bound for prime pairs and the special case $k=2$ of Gallagher's Theorem. Since the prime pair $p' = p+k$ corresponds not only to the prime 2-tuple $(n,n+k)$ but any shifted tuple $(n+j, n+j+k)$, we have 
\[ \sum_{\substack{N<p,p'\le 2N\\ 0<p'-p\le h}} 1 < \frac1h \sum_{\substack{|\mathcal{H}|=2\\ \mathcal{H} \subset \{1, 2, \ldots , \lfloor 2h\rfloor\}}}\sum_{\substack{\frac N2 < n < 3N \\ (\mathcal{P}_{\mathcal H}(n), P(N^{\frac12}))=1}} 1,\]
which by \eqref{eq:5.9} is
\[ \ll \frac {N}{h(\log N)^2} \sum_{\substack{|\mathcal{H}|=2\\ \mathcal{H} \subset \{1, 2, \ldots , \lfloor 2h\rfloor\}}} \mathfrak{S}(\mathcal H)  \quad \ll  h \frac{N}{(\log N)^2},\]
by \eqref{eq:2.9}, which is equivalent to Theorem 2.

\vspace{1cm} \footnotesize D. A. Goldston  \,\,\,
(goldston@math.sjsu.edu)

Department of Mathematics

San Jose State University

San Jose, CA 95192

 USA \\

J. Pintz \,\,\, (pintz@renyi.hu)

R\'enyi Mathematical Institute of the Hungarian Academy of Sciences

H-1364 Budapest

P.O.B. 127

Hungary \\

C. Y. Y{\i}ld{\i}r{\i}m \,\,\, (yalciny@boun.edu.tr)

Department of Mathematics

Bo\~{g}azi\c{c}i University

Bebek, Istanbul 34342

Turkey \\

\end{document}